\documentclass[11pt, twoside]{article}
\usepackage{latexsym}
\usepackage{amsmath}
\usepackage{amssymb}
\usepackage[all]{xy}
\usepackage{amsfonts}
\usepackage{verbatim}
\usepackage{amsthm}
\usepackage{bm}
\usepackage{mathrsfs}
\usepackage{epsfig}
\usepackage{xy}
\usepackage{array}
\usepackage{stmaryrd}
\usepackage{graphicx,color}
\usepackage{xcolor}
\usepackage{tikz}
\usetikzlibrary{arrows,calc}
\usepackage{etex}
\usepackage{mathdots}
\usepackage{float}
\usepackage{graphics}
\usepackage{pdflscape}
\usepackage{CJK}
\usepackage{anysize,hyperref}
\input xypic
\xyoption{all}
\usepackage{extarrows}
\usepackage[perpage,symbol]{footmisc}
\usepackage{setspace}
\renewcommand{\cal}{\mathcal}

\def\F{\mathfrak{F}}

\def\P{\mathcal{P}}
\def\I{\mathcal{I}}
\def\E{\mathbb{E}}
\def\B{\mathcal{B}}
\def\C{\mathcal{C}}
\def\s{\mathfrak{s}}

\def\op{^\mathrm{op}}
\def\Ab{\mathit{Ab}}
\def\del{\delta}
\def\dr{\ar@{->}[r]}

\def\X{\mathcal{X}}
\def\Y{\mathcal{Y}}
\def\Z{\mathcal{Z}}

\def\Ext{\mbox{Ext}}
\def\Hom{\mbox{Hom}}

\begin{document}
\baselineskip=15pt
\title{\Large{\bf Filtered objects in extriangulated categories$^\bigstar$\footnotetext{\hspace{-1em}$^\bigstar$Panyue Zhou is supported by the National Natural Science Foundation of China (Grants No. 11901190 and 11671221) and the Hunan Provincial Natural Science Foundation of China (Grants No. 2018JJ3205).}}}
\medskip
\author{Panyue Zhou}

\date{}

\maketitle
\def\blue{\color{blue}}
\def\red{\color{red}}

\newtheorem{theorem}{Theorem}[section]
\newtheorem{lemma}[theorem]{Lemma}
\newtheorem{corollary}[theorem]{Corollary}
\newtheorem{proposition}[theorem]{Proposition}
\newtheorem{conjecture}{Conjecture}
\theoremstyle{definition}
\newtheorem{definition}[theorem]{Definition}
\newtheorem{question}[theorem]{Question}
\newtheorem{remark}[theorem]{Remark}
\newtheorem{remark*}[]{Remark}
\newtheorem{example}[theorem]{Example}
\newtheorem{example*}[]{Example}
\newtheorem{condition}[theorem]{Condition}
\newtheorem{condition*}[]{Condition}
\newtheorem{construction}[theorem]{Construction}
\newtheorem{construction*}[]{Construction}

\newtheorem{assumption}[theorem]{Assumption}
\newtheorem{assumption*}[]{Assumption}

\baselineskip=17pt
\parindent=0.5cm

\begin{abstract}
\baselineskip=16pt
Let $R$ be an artin ring and $\Theta=\{\Theta(1),\Theta(2),\cdots,\Theta(n)\}$ be
a family of objects in an artin extriangulated $R$-category $(\C,\E,\s)$ such that
$\E(\Theta(j),\Theta(i))=0$ for all $j\geq i$. In this paper, we show that
the class $\P(\Theta)$ of the $\Theta$-projective objects is a precovering class and the class $\I(\Theta)$ of the $\Theta$-injective objects is a preenveloping one in $\C$.
Furthermore, if $\C$ has enough projectives and enough injectives, we show that the subcategory $\F(\Theta)$ of $\Theta$-filtered objects is functorially finite in $\C$. As an appliacation, this generalizes the works by Ringel in a module category case and Mendoza-Santiago in a triangulated category case.\\[0.5cm]
\medskip
\end{abstract}

\pagestyle{myheadings}
\markboth{\rightline {\scriptsize   Panyue Zhou}}
         {\leftline{\scriptsize  Filtered objects in extriangulated categories}}

\section{Introduction}
Let $\X$ and $\Y$ be classes of objects in an additive category $\C.$ A morphism $f:X\to C$ in $\C$ is said to be an
\emph{$\X$-precover} of $C$ if $X\in\X$ and $\Hom_\C(X',f)\colon\Hom_\C(X',X)\to\Hom_\C(X',C)$ is surjective, for any $ X'\in\X.$
 If any $C\in\Y$ admits an $\X$-precover, then $\X$ is called a \emph{precovering} class in $\C.$ By duality, we can define the notions of an \emph{$\X$-preenvelope} of $C$ and a \emph{preenveloping} class in $\C.$ It is said that $\X$ is \emph{functorially finite} in $\C$ if $\X$ is both precovering and preenveloping in $\C.$
For more details, see \cite{AR}.

Let $\Lambda$ be an Artin algebra and ${\rm mod}(\Lambda)$ the category of finitely generated left $\Lambda$-modules.
Assume that $\Theta=\{\Theta(1),\Theta(2),\cdots,\Theta(n)\}$ is a finite set of $\Lambda$-modules with
${\rm Ext}^1_{\Lambda}(\Theta(j), \Theta(i))=0$ for $j\geq i$. We denote by $\F(\Theta)$ the full subcategory of ${\rm mod}(\Lambda)$ of modules having a filtration with factors in $\Theta$. Thus, $M$ belongs to $\F(\Theta)$ if and only if $M$ has
submodules $$0=M_0\subseteq M_1\subseteq\cdots\subseteq M_t=M$$ such that $M_s/M_{s-1}$ is isomorphic to
a module in $\Theta$.

Ringel showed the following well-known result. Note that in this way we obtain a large variety of
functorially finite subcategories of ${\rm mod}(\Lambda)$ which usually will not be closed under
submodules or factor modules.
\begin{theorem}\emph{\cite[Theorem 1]{R}}
The subcategory $\F(\Theta)$ is functorially finite in ${\rm mod}(\Lambda)$.
\end{theorem}

Let $R$ be an artin ring and $\X$ be a class of objects in an artin triangulated $R$-category $\C$.
It is said that an object $M\in\C$ admits an
\emph{$\X$-filtration} if there exists a family of
distinguished triangles $$\eta=\{\eta_{i}\colon M_{i-1}\to M_{i}\to X_{i}\to M_{i-1}[1]\}_{i=0}^{n}$$ such that $M_{-1}=0=X_{0}$, $M_{n}=M$ and
$X_{i}\in \X$ for any $i\geq 1$.
We denote by  $\F(\X)$ the class of objects
$M\in \C$ for which there exists an $\X$-filtration.

Mendoza-Santiago proved a triangulated version of Ringel's theorem.

\begin{theorem}\emph{\cite[Theorem 4.10]{MS}}
 Let $\Theta=\{\Theta(1),\Theta(2),\cdots,\Theta(n)\}$ be
a family of objects in an artin triangulated $R$-category $\C$ such that
${\rm Hom}_{\C}(\Theta(j),\Theta(i)[1])=0$ for all $j\geq i$. Then
the subcategory $\F(\Theta)$ is functorially finite in $\C$.
\end{theorem}

Extriangulated categories were recently introduced by Nakaoka and Palu \cite{NP} by
extracting those properties of $\Ext^1$ on exact categories (which is itself a generalisation of the concept
of a module category and an abelian category) and on triangulated categories
that seem relevant from the point of view of cotorsion pairs. In particular, exact
categories and triangulated categories are extriangulated categories. There are a lot of examples of extriangulated categories which are neither exact triangulated categories nor triangulated categories, see \cite{NP,ZZ,HZZ}.
Hence, many results hold on exact categories and triangulated categories can be unified in the same framework.

Motivated
by this idea, we replace $\E$-triangles in extriangulated categories with distinguished triangles in triangulated categories to define the subcategory $\F(\Theta)$. We can unify Ringel's result and Mendoza-Santiago's result under the framework of extriangulated categories.
Our first main result is the following.
\begin{theorem}\emph{(See Theorem \ref{main1} for more details)}
Let $(\C,\E,\s)$ be an artin extriangulated $R$-category with enough projectives and enough injectives,
and $\Theta=\{\Theta(1),\Theta(2),\cdots,\Theta(n)\}$ be
a family of objects in $\C$ such that
$\E(\Theta(j),\Theta(i))=0$ for all $j\geq i$. Then
the subcategory $\F(\Theta)$ is functorially finite in $\C$.
\end{theorem}

Moreover, we also define the notions of $\Theta$-projective objects and $\Theta$-injective objects in an extriangulated category $\C$, see Definition \ref{def1}.  Our second main result is the following.
\begin{theorem}\emph{(See Theorem \ref{main2} for more details)}
Let be $\Theta=\{\Theta(1),\Theta(2),\cdots,\Theta(n)\}$ be
a family of objects in an artin extriangulated $R$-category $(\C,\E,\s)$ such that
$\E(\Theta(j),\Theta(i))=0$ for all $j\geq i$. Then
the class $\P(\Theta)$ of the $\Theta$-projective objects is a precovering class and the class $\I(\Theta)$ of the $\Theta$-injective objects is a preenveloping one in $\C$.
\end{theorem}
This paper is organised as follows:

In Section 2, we review some elementary definitions on extriangulated categories.

In Section 3, we show the first main result and second one.

\section{Preliminaries}
Let us briefly recall the definition and basic properties of extriangulated categories from \cite{NP}.

Let $\C$ be an additive category. Suppose that $\C$ is equipped with a biadditive functor $$\E\colon\C\op\times\C\to\Ab,$$
where $\Ab$ is the category of abelian groups. For any pair of objects $A,C\in\C$, an element $\delta\in\E(C,A)$ is called an {\it $\E$-extension}. Thus formally, an $\E$-extension is a triplet $(A,\delta,C)$.
Let $(A,\del,C)$ be an $\E$-extension. Since $\E$ is a bifunctor, for any $a\in\C(A,A')$ and $c\in\C(C',C)$, we have $\E$-extensions
$$ \E(C,a)(\del)\in\E(C,A')\ \ \text{and}\ \ \ \E(c,A)(\del)\in\E(C',A). $$
We abbreviately denote them by $a_\ast\del$ and $c^\ast\del$.
For any $A,C\in\C$, the zero element $0\in\E(C,A)$ is called the \emph{spilt $\E$-extension}.

\begin{definition}{\cite[Definition 2.3]{NP}}
Let $(A,\del,C),(A',\del',C')$ be any pair of $\E$-extensions. A {\it morphism} $$(a,c)\colon(A,\del,C)\to(A',\del',C')$$ of $\E$-extensions is a pair of morphisms $a\in\C(A,A')$ and $c\in\C(C,C')$ in $\C$, satisfying the equality
$ a_\ast\del=c^\ast\del'. $
Simply we denote it as $(a,c)\colon\del\to\del'$.
\end{definition}

\begin{definition}{\cite[Definition 2.6]{NP}}
Let $\delta=(A,\delta,C),\delta^{\prime}=(A^{\prime},\delta^{\prime},C^{\prime})$ be any pair of $\mathbb{E}$-extensions. Let
\[ C\overset{\iota_C}{\longrightarrow}C\oplus C^{\prime}\overset{\iota_{C^{\prime}}}{\longleftarrow}C^{\prime} \]
and
\[ A\overset{p_A}{\longleftarrow}A\oplus A^{\prime}\overset{p_{A^{\prime}}}{\longrightarrow}A^{\prime} \]
be coproduct and product in $\B$, respectively. Since $\mathbb{E}$ is biadditive, we have a natural isomorphism
\[ \mathbb{E}(C\oplus C^{\prime},A\oplus A^{\prime})\cong \mathbb{E}(C,A)\oplus\mathbb{E}(C,A^{\prime})\oplus\mathbb{E}(C^{\prime},A)\oplus\mathbb{E}(C^{\prime},A^{\prime}). \]

Let $\delta\oplus\delta^{\prime}\in\mathbb{E}(C\oplus C^{\prime},A\oplus A^{\prime})$ be the element corresponding to $(\delta,0,0,\delta^{\prime})$ through the above isomorphism. This is the unique element which satisfies
\begin{eqnarray*}
\mathbb{E}(\iota_C,p_A)(\delta\oplus\delta^{\prime})=\delta&,&\mathbb{E}(\iota_C,p_{A^{\prime}})(\delta\oplus\delta^{\prime})=0,\\
\mathbb{E}(\iota_{C^{\prime}},p_A)(\delta\oplus\delta^{\prime})=0&,&\mathbb{E}(\iota_{C^{\prime}},p_{A^{\prime}})(\delta\oplus\delta^{\prime})=\delta^{\prime}.
\end{eqnarray*}
\end{definition}

\medskip

Let $A,C\in\C$ be any pair of objects. Sequences of morphisms in $\C$
$$\xymatrix@C=0.7cm{A\ar[r]^{x} & B \ar[r]^{y} & C}\ \ \text{and}\ \ \ \xymatrix@C=0.7cm{A\ar[r]^{x'} & B' \ar[r]^{y'} & C}$$
are said to be {\it equivalent} if there exists an isomorphism $b\in\C(B,B')$ which makes the following diagram commutative.
$$\xymatrix{
A \ar[r]^x \ar@{=}[d] & B\ar[r]^y \ar[d]_{\simeq}^{b} & C\ar@{=}[d]&\\
A\ar[r]^{x'} & B' \ar[r]^{y'} & C &}$$

We denote the equivalence class of $\xymatrix@C=0.7cm{A\ar[r]^{x} & B \ar[r]^{y} & C}$ by $[\xymatrix@C=0.7cm{A\ar[r]^{x} & B \ar[r]^{y} & C}]$.

For any $A,C\in\C$, we denote $0=[
A \xrightarrow{\binom{1}{0}} A\oplus C \xrightarrow{(0,1)} C ]$.

For any two equivalence classes, we denote
$$ [A\overset{x}{\longrightarrow}B\overset{y}{\longrightarrow}C]\oplus [A^{\prime}\overset{x^{\prime}}{\longrightarrow}B^{\prime}\overset{y^{\prime}}{\longrightarrow}C^{\prime}]=[A\oplus A^{\prime}\overset{x\oplus x^{\prime}}{\longrightarrow}B\oplus B^{\prime}\overset{y\oplus y^{\prime}}{\longrightarrow}C\oplus C^{\prime}]. $$

\begin{definition}{\cite[Definition 2.9]{NP}}
Let $\s$ be a correspondence which associates an equivalence class $\s(\del)=[\xymatrix@C=0.7cm{A\ar[r]^{x} & B \ar[r]^{y} & C}]$ to any $\E$-extension $\del\in\E(C,A)$. This $\s$ is called a {\it realization} of $\E$, if it satisfies the following condition:
\begin{itemize}
\item Let $\del\in\E(C,A)$ and $\del'\in\E(C',A')$ be any pair of $\E$-extensions, with $$\s(\del)=[\xymatrix@C=0.7cm{A\ar[r]^{x} & B \ar[r]^{y} & C}],\ \ \ \s(\del')=[\xymatrix@C=0.7cm{A'\ar[r]^{x'} & B'\ar[r]^{y'} & C'}].$$
Then, for any morphism $(a,c)\colon\del\to\del'$, there exists $b\in\C(B,B')$ which makes the following diagram commutative.
$$\xymatrix{
A \ar[r]^x \ar[d]^a & B\ar[r]^y \ar[d]^{b} & C\ar[d]^c&\\
A'\ar[r]^{x'} & B' \ar[r]^{y'} & C' &}$$

\end{itemize}
In the above situation, we say that the triplet $(a,b,c)$ realizes $(a,c)$.
\end{definition}

\begin{definition}{\cite[Definition 2.10]{NP}}
A realization $\s$ of $\E$ is called \emph{additive} if it satisfies the following conditions.
\begin{itemize}
\item[(1)] For any $A,C\in\C$, the split $\E$-extension $0\in\E(C,A)$ satisfies $\s(0)=0$.
\item[(2)] For any pair of $\E$-extensions $\delta\in\E(C,A)$ and $\delta'\in\E(C',A')$,
$$\s(\delta\oplus\delta')=\s(\delta)\oplus\s(\delta')$$
holds.
\end{itemize}
\end{definition}

\begin{definition}{\cite[Definition 2.12]{NP}}
A triplet $(\C,\E,\s)$ is called an \emph{externally triangulated category} (or \emph{extriangulated category} for short) if it satisfies the following conditions:
\begin{itemize}
\item[{\rm (ET1)}] $\E\colon\C\op\times\C\to\Ab$ is a biadditive functor.
\item[{\rm (ET2)}] $\s$ is an additive realization of $\E$.
\item[{\rm (ET3)}] Let $\del\in\E(C,A)$ and $\del'\in\E(C',A')$ be any pair of $\E$-extensions, realized as
$$ \s(\del)=[\xymatrix@C=0.7cm{A\ar[r]^{x} & B \ar[r]^{y} & C}],\ \ \s(\del')=[\xymatrix@C=0.7cm{A'\ar[r]^{x'} & B' \ar[r]^{y'} & C'}]. $$
For any commutative square
$$\xymatrix{
A \ar[r]^x \ar[d]^a & B\ar[r]^y \ar[d]^{b} & C&\\
A'\ar[r]^{x'} & B' \ar[r]^{y'} & C' &}$$
in $\C$, there exists a morphism $(a,c)\colon\del\to\del'$ satisfying $cy=y'b$.
\item[{\rm (ET3)$\op$}] Dual of (ET3).
\item[{\rm (ET4)}] Let $(A,\del,D)$ and $(B,\del',F)$ be $\E$-extensions realized by
$$\xymatrix@C=0.7cm{A\ar[r]^{f} & B \ar[r]^{f'} & D}\ \ \text{and}\ \ \ \xymatrix@C=0.7cm{B\ar[r]^{g} & C \ar[r]^{g'} & F}$$
respectively. Then there exist an object $E\in\C$, a commutative diagram
$$\xymatrix{A\ar[r]^{f}\ar@{=}[d]&B\ar[r]^{f'}\ar[d]^{g}&D\ar[d]^{d}\\
A\ar[r]^{h}&C\ar[d]^{g'}\ar[r]^{h'}&E\ar[d]^{e}\\
&F\ar@{=}[r]&F}$$
in $\C$, and an $\E$-extension $\del^{''}\in\E(E,A)$ realized by $\xymatrix@C=0.7cm{A\ar[r]^{h} & C \ar[r]^{h'} & E},$ which satisfy the following compatibilities.
\begin{itemize}
\item[{\rm (i)}] $\xymatrix@C=0.7cm{D\ar[r]^{d} & E \ar[r]^{e} & F}$  realizes $f'_{\ast}\del'$,
\item[{\rm (ii)}] $d^\ast\del''=\del$,

\item[{\rm (iii)}] $f_{\ast}\del''=e^{\ast}\del'$.
\end{itemize}

\item[{\rm (ET4)$\op$}]  Dual of (ET4).
\end{itemize}
\end{definition}

We use the following terminology.
\begin{definition}{\cite{NP}}
Let $(\C,\E,\s)$ be an extriangulated category.
\begin{itemize}
\item[(1)] A sequence $A\xrightarrow{~x~}B\xrightarrow{~y~}C$ is called a {\it conflation} if it realizes some $\E$-extension $\del\in\E(C,A)$.
    In this case, $x$ is called an {\it inflation} and $y$ is called a {\it deflation}.

\item[(2)] If a conflation  $A\xrightarrow{~x~}B\xrightarrow{~y~}C$ realizes $\delta\in\mathbb{E}(C,A)$, we call the pair $( A\xrightarrow{~x~}B\xrightarrow{~y~}C,\delta)$ an {\it $\E$-triangle}, and write it in the following way.
$$A\overset{x}{\longrightarrow}B\overset{y}{\longrightarrow}C\overset{\delta}{\dashrightarrow}$$
We usually do not write this $``\delta"$ if it is not used in the argument.

\item[(3)] Let $A\overset{x}{\longrightarrow}B\overset{y}{\longrightarrow}C\overset{\delta}{\dashrightarrow}$ and $A^{\prime}\overset{x^{\prime}}{\longrightarrow}B^{\prime}\overset{y^{\prime}}{\longrightarrow}C^{\prime}\overset{\delta^{\prime}}{\dashrightarrow}$ be any pair of $\E$-triangles. If a triplet $(a,b,c)$ realizes $(a,c)\colon\delta\to\delta^{\prime}$, then we write it as
$$\xymatrix{
A \ar[r]^x \ar[d]^a & B\ar[r]^y \ar[d]^{b} & C\ar@{-->}[r]^{\del}\ar[d]^c&\\
A'\ar[r]^{x'} & B' \ar[r]^{y'} & C'\ar@{-->}[r]^{\del'} &}$$
and call $(a,b,c)$ a {\it morphism of $\E$-triangles}.

\item[(4)] An object $P\in\C$ is called {\it projective} if
for any $\E$-triangle $A\overset{x}{\longrightarrow}B\overset{y}{\longrightarrow}C\overset{\delta}{\dashrightarrow}$ and any morphism $c\in\C(P,C)$, there exists $b\in\C(P,B)$ satisfying $yb=c$.
We denote the subcategory of projective objects by $\cal P\subseteq\C$. Dually, the subcategory of injective objects is denoted by $\cal I\subseteq\C$.

\item[(5)] We say that $\C$ {\it has enough projective objects} if
for any object $C\in\C$, there exists an $\E$-triangle
$A\overset{x}{\longrightarrow}P\overset{y}{\longrightarrow}C\overset{\delta}{\dashrightarrow}$
satisfying $P\in\cal P$. We can define the notion of having enough injectives dually.

\end{itemize}
\end{definition}

Now we give some basic facts that will be used later.

\begin{lemma}\label{y0}\emph{\cite[Corollary 3.5]{NP}}
 Assume that $\C$ is an extriangulated category.
Let
$$\xymatrix{A\ar[d]^{a}\ar[r]^{x} & B\ar[r]^y\ar[d]^{b} &
C\ar[d]^c\ar@{-->}[r]^{\delta}&\\
 A'\ar[r]^{x'} & B'\ar[r]^{y'} & C'\ar@{-->}[r]^{\delta'}&}$$
Then $a$ factors through $x$ if and only if $a_{\ast}\del=c^{\ast}\del'$ if and only if $c$ factors through $y'$.
In particular, we obtain
$x$ is a split monomorphism $\Longleftrightarrow$ $\delta$ splits $\Longleftrightarrow$ $y$ is a split epimorphism.
\end{lemma}

\begin{lemma}\label{y2}\emph{\cite[Proposition  1.20]{LN}}
Let $\C$ be an extriangulated catgeory and
$$A\overset{f}{\longrightarrow}B\overset{g}{\longrightarrow}C\overset{}{\dashrightarrow}$$
be an $\E$-triangle in $\C$. Then for  every morphism $a\colon A\to X$, there exists a commutative diagram
$$\xymatrix{A\ar[r]^f \ar[d]_{a} & B\ar[r]^g \ar[d]^{b} & C\ar@{-->}[r] \ar@{=}[d]&\\
X\ar[r]^{u} &Y \ar[r]^{v} &C \ar@{-->}[r]&}$$
of $\E$-triangles.
\end{lemma}

Assume that $(\C, \E, \s)$ is an extriangulated category. By Yoneda's lemma, any $\E$-extension
$\del\in\E(C,A)$ induces natural transformations
$$\del_{\sharp}\colon\C(-,C)\Rightarrow\E(-,A)\ \ \textrm{and} \ \ \del^{\sharp}\colon \C(A,-)\Rightarrow\E(C,-).$$
For any $X\in\C$, these $(\del_{\sharp})_X $ and $\del_{X}^{\sharp}$ are given as follows:

(1) $(\del_{\sharp})_X\colon\C(X,C)\to \E(X,A); f\mapsto f^{\ast}\del$.

(2) $\del_{X}^{\sharp}\colon \C(A,X)\to \E(C,X); g\mapsto g_{\ast}\del$.

When there is no danger of confusion, we will sometimes instead of
$$\Hom_\C(X,A)\xrightarrow{~\Hom_\C(X,\ f)~}\Hom_\C(X,B)$$
 write the simplified form:
$$\C(X,A)\xrightarrow{~\C(X,\ f)~}\C(X,B).$$

\begin{lemma}\label{a0}
Let $\C$ be an extriangulated category, $$\xymatrix{A\ar[r]^{x}&B\ar[r]^{y}&C\ar@{-->}[r]^{\delta}&}$$
an $\E$-triangle. Then we have the following long exact sequence:
$$\C(-, A)\xrightarrow{\C(-,x)}\C(-, B)\xrightarrow{\C(-,y)}\C(-, C)\xrightarrow{\delta^{\sharp}_-}
\E(-, A)\xrightarrow{\E(-,x)}\E(-, B)\xrightarrow{\E(-,y)}\E(-, C);
$$
$$\C(C,-)\xrightarrow{\C(y,-)}\C(B,-)\xrightarrow{\C(x,-)}\C(A,-)\xrightarrow{\delta_{\sharp}^-}
\E(C,-)\xrightarrow{\E(y,-)}\E(B,-)\xrightarrow{\E(x,-)}\E(A,-)
.$$
\end{lemma}

\proof This follows from Proposition 3.3 and Proposition 3.11 in \cite{NP}. \qed
\bigskip

Let $R$ be a commutative ring. Recall that an \emph{$R$-category} is a category $\C$ satisfying
the following two conditions: (a) for any pair $X,Y$ of objects in $\C,$ the set of morphisms $\Hom_\C(X,Y)$ is an $R$-module, and (b) the composition of morphisms in $\C$ is $R$-bilinear.
 An $R$-category $\C$ is called \emph{Hom-finite} if $\Hom_\C(X,Y)$ is a finitely generated $R$-module, for any $X,Y\in \C.$
An extriangulated $R$-category $(\C,\E,\s)$ is called \emph{Ext-finite} if $\E(X,Y)$ is a finitely generated $R$-module, for any $X,Y\in \C.$ An additive category is \emph{Krull-Schmidt} if each of its objects
is the direct sum of finitely many objects with local endomorphism rings. It
follows that these finitely many objects are indecomposable and determined up to
isomorphism.

\begin{definition}
An artin extriangulated $R$-category is an $R$-category for some artin
ring $R$, which is a Hom-finite, Ext-finite,  Krull-Schmidt and extriangulated category.
\end{definition}

Now we give some examples of artin extriangulated $R$-categories.
\begin{example}Let $R$ be an artin ring.
\begin{itemize}
\item[\rm(1)] Let $\Lambda$ be an artin $R$-algebra. It is well known that ${\rm mod}(\Lambda)$ is an artin exact
$R$-category. It can be viewed as an artin extriangulated $R$-category.

\item[\rm(2)]  Let $\Lambda$ be an artin $R$-algebra. It is well-known
that the bounded derived category $D^b(\Lambda)$ of complexes in ${\rm mod}(\Lambda),$ is an
artin triangulated $R$-category, see \cite[Theorem B.2]{CYZ}.
It can be viewed as an artin extriangulated $R$-category.
We know that the extension closed subcategory of a triangulated category (See \cite[Remark 2.18]{NP}) is an extriangulated category, but it is not exact nor triangulated.
Thus we take an extension closed subcatgory $\B$ of $D^b(\Lambda)$, we get that $\B$ is an
artin extriangulated $R$-category.
\end{itemize}
\end{example}

The following construction is usually called the \emph{universal extension}.
\begin{lemma}\label{y1} Let $\C$ be an artin extriangulated $R$-category,
 and let $A,C\in \C$ be such that
$\E(C, A)\neq 0$. Then we have the following conditions hold.
\begin{itemize}
\item [\rm (a)] There exists a not splitting $\E$-triangle in
$\C$
$$\eta_{C,A}\colon \xymatrix{A^{n}\ar[r]^{f} &B\ar[r]^{g} & C\ar@{-->}[r]&}$$ such that
$\delta^{\sharp}_A\colon{\rm Hom}_{\C}(A^{n},A) \to \E(C,A)$ is surjective, where
$n:=\ell_{R}(\E(C,A))$ is finite $R$-length of $\E(C,A)$.

\item [\rm (b)] If $\E(A, A)=0$ then $\E(B, A)=0$.
\end{itemize}
\end{lemma}

\proof  The proof given in \cite[Lemma 3.4]{MS} can be adapted to the context of extriangulated
categories.   \qed

\section{Filtered objects}
Let $(\C,\E,\s)$ be an extriangulated category and $\X$ be a class of objects in
$\C$. It is said that an object $M\in\C$ admits an
\emph{$\X$-filtration} if there exists a family of
$\E$-triangles $$\eta=\{\eta_{i}\colon M_{i-1}\to M_{i}\to X_{i}\dashrightarrow\}_{i=0}^{n}$$ such that $M_{-1}=0=X_{0}$, $M_{n}=M$ and
$X_{i}\in \X$ for $i\geq 1$. In this case,  we define two
lengths: $\ell_{\X,\eta}(M):=n$ and
$\ell_{\X}(M):=\min\{\ell_{\X,\eta}(M)\mid \eta
\,\,\text{is an}\,\,\X\text{-filtration of}\, \,M\}$.
We denote by  $\F(\X)$ the class of objects
$M\in \C$ for which there exists an $\X$-filtration.

We have two left and two right perpendicular classes of $\X$ to be
$$^\perp\X:=\{A\in\C~|~\Hom_{\C}(A,X)=0,~\textrm{for any}~ X\in\X\};$$
$$\X^\perp:=\{A\in\C~|~\Hom_{\C}(X,A)=0,~\textrm{for any}~ X\in\X\};$$
$$^{_{\E}\perp}\X:=\{A\in\C~|~\E(A,X)=0,~\textrm{for any}~ X\in\X\}.$$
$$\X^{_{\E}\perp}:=\{A\in\C~|~\E(X,A)=0,~\textrm{for any}~ X\in\X\}.$$
For two classes $\X$ and $\Y$ of objects in $\C$, we denote by $\X\ast\Y$ the class of objects $A\in\C$ for which
there exists an $\E$-triangle
$X\to A\to Y\dashrightarrow$
in $\C$ with $X\in\X$ and $Y\in\Y$. Moreover, it is said that $\X$ is closed under extensions if $\X\ast\X\subseteq\X$.

\begin{remark}\label{rem1}
Let $(\C,\E,\s)$ be an extriangulated category and $\X$ a class of object in $\C$. Then we have the following statements hold.
\begin{enumerate}
\item [(1)] $\F(\X)=\bigcup\limits_{n\in
\mathbb{N}}\F_{n}(\X)$, where
$\F_{0}(\X):=\{0\}$ and
$\F_{n}(\X):=\F_{n-1}(\X)\ast \X$  for
 any $n\geq 1$.
\item [(2)] $\ell_{\X}(M)=\min\,\{n\in \mathbb{N}\mid M\in
\F_{n}(\X)\}$ for any $M\in \F(\X).$
\end{enumerate}
\end{remark}

\begin{lemma}\label{lem2}
Let $\C$ be an extriangulated category and $\X$ a class of objects in $\C$.
Then the class $\F(\X)$ is closed under extensions. In particular, it is closed under finite direct sums.
\end{lemma}

\proof Assume that  $A\to B\to C\dashrightarrow$ is an $\E$-triangle in $\C$ with $A$ and $C$ in $\F(\X)$.
The proof will be completed by induction on  $n:=\ell_{\X}(C)$.

If $C=0$, we have 
$A\simeq B$ and then $B\in\F(\X)$.

If $\ell_{\X}(C)=1$, we have $C\simeq X\in \X,$ and therefore
an $\X$-filtration of $B$ can be done by adding the $\E$-triangle $A\to B\to C\dashrightarrow$ to an
$\X$-filtration of $A$.

Now we suppose $\ell_{\X}(C)>1$. We consider a minimal $\X$-filtration of $C$,
$$\{\eta_{i}\colon\xymatrix{C_{i-1}\ar[r] &
C_{i}\ar[r] & X_{i}\ar@{-->}[r]&}\}_{i=0}^{n}.$$  By (ET4)$\op$,
we have the following commutative diagram 
$$\xymatrix{A\ar[r]\ar@{=}[d] & B_{n-1}\ar[r]\ar[d] & C_{n-1}\ar@{-->}[r]\ar[d]&\\
A\ar[r] & B\ar[r]\ar[d] & C\ar@{-->}[r]\ar[d] &\\
& X_{n}\ar@{-->}[d]\ar@{=}[r] & X_{n}\ar@{-->}[d]\\&&}$$
of $\E$-triangles. It follows that $\ell_{\X}(C_{n-1})<\ell_{\X}(C)$. Applying induction to the first row of
the above diagram, we have that $B_{n-1}\in\F(\X)$.
Hence, an $\X$-filtration of $B$ is given by adding the $\E$-triangle
$B_{n-1}\to B\to X_{n}\dashrightarrow$ to an $\X$-filtration of $B_{n-1}$.
This completes the proof.\qed

\begin{lemma}\label{lem3}
Let $\C$ be an extriangulated category, $\Y$ and $\Z$ two classes of objects in $\C$.
\begin{itemize}
\item[\rm (1)] If ${\rm Hom}_{\C}(\Y,\Z)=0$, then
${\rm Hom}_{\C}(\F(\Y),\F(\Z))=0$.
\item[\rm (2)] If $\E(\Y,\Z)=0$, then
$\E(\F(\Y),\F(\Z))=0$.
\end{itemize}
\end{lemma}

\proof (1)
Suppose that $N\in\F(\Y)$ and $M\in\F(\Z)$. We will show that $\Hom_{\C}(N,M)=0$ by induction on
$\ell_{\Y}(N)$. Without loss of generality, we also can assume  $M\neq0$ and $N\neq 0$.

If $\ell_{\Y}(N)=1$ then $N\simeq Y\in \Y$ and so it can be seen that
$\Hom_{\C}(N,M)=0$ by induction on $\ell_{\Z}(M)$.

Suppose that $n:=\ell_{\Y}(N)>1$. Then there exists an $\E$-triangle
$$\eta_{n}\colon~\xymatrix{N_{n-1}\ar[r] & N\ar[r] &
Y_{n}\ar@{-->}[r] &}$$ such that $N_{n-1}\in
\F(\Y)$, $Y_{n}\in \Y$ and
$\ell_{\Y}(N_{n-1})=n-1$. Applying the functor $\Hom_{\C}(-,M)$ to
the above $\E$-triangle $\eta_{n}$, we have the following exact sequence
$$\xymatrix{\Hom_{\C}(Y_{n},M)\ar[r] &
\Hom_{\C}(N,M)\ar[r] & \Hom_{\C}(N_{n-1},M).}$$ By
induction, we get that $\Hom_{\C}(N_{n-1},M)=0=\Hom_{\C}(Y_{n},M),$ and so
$\Hom_{\C}(N,M)=0$.

(2) The proof is similar to (1).  \qed

\begin{lemma} \label{lem4}
Let $\C$ be an extriangulated category and $\X$ a class of objects in $\C$. Then
$$^{\perp}\X=\,^{\perp}\F(\X)~~\textrm{and}~~^{_{\E}\perp}\X=\,^{_{\E}\perp}\F(\X).$$
\end{lemma}

\proof We only show that $^{\perp}\X=\,^{\perp}\F(\X)$. Similarly, one can show that
$^{_{\E}\perp}\X=\,^{_{\E}\perp}\F(\X)$.

It suffices to show that ${}^{\perp} \X\subseteq{}^{\perp}\F(\X),$ since the other side inclusion ${}^{\perp}\F(\X)
\subseteq {}^{\perp} \X$ follows easily from the fact that $\X\subseteq \F(\X).$

Suppose $M\in{^{\perp}\X}$ and  $N\in\F(\X)$. By Lemma \ref{lem3}, we have
that $\Hom_{\C}(M,N)=0$ since $\Hom_{\C}(M,\X)=0$. It follows that $Y\in{}^{\perp}\F(\X)$
and then $^{\perp}\X\subseteq\,^{\perp}\F(\X)$ as required.  \qed

\begin{lemma}\label{lem5}
Let $\C$ be an extriangulated category.
If there are two $\E$-triangles
$$Z\to Y\to U\dashrightarrow~~\textrm{and}~~Y\to X\to V\dashrightarrow$$ such that
$\E(V,U)=0$, then there exist two $\E$-triangles as follows
$$Z\to W\to V\dashrightarrow~~\textrm{and}~~W\to X\to U\dashrightarrow.$$
\end{lemma}

\proof
By (ET4), we have the following commutative diagram
$$\xymatrix{Z\ar[r]\ar@{=}[d] & Y\ar[r]\ar[d] & U\ar@{-->}[r]\ar[d]^{u}&\\
Z\ar[r] & X\ar[r]\ar[d] & C\ar@{-->}[r]\ar[d]^{v} &\\
& V\ar@{-->}[d]\ar@{=}[r] & V\ar@{-->}[d]\\&&}$$
of $\E$-triangles. Using the fact that $\E(V,U)=0,$ it follows that $\eta\colon U\xrightarrow{~u~}C\xrightarrow{~v~} V\dashrightarrow$ splits, and hence we get the following splitting $\E$-triangle $\eta'\colon V\xrightarrow{~u'~} C\xrightarrow{~v'~}U\dashrightarrow$. By (ET4)$\op$, we
obtain the following commutative diagram
$$\xymatrix{Z\ar[r]\ar@{=}[d] & W\ar[r]\ar[d] &V\ar@{-->}[r]\ar[d]^{u'}&\\
Z\ar[r] & X\ar[r]\ar[d] & C\ar@{-->}[r]\ar[d]^{v'} &\\
& U\ar@{-->}[d]\ar@{=}[r] & U\ar@{-->}[d]\\&&}$$
of $\E$-triangles. Thus the required $\E$-triangles are $Z\to W\to V\dashrightarrow~~\textrm{and}~~W\to X\to U\dashrightarrow.$ \qed

\begin{lemma}\label{lem6}
Let $\C$ be an extriangulated category and $A\in\C$ with $\E(A,A)=0$, and let
$$\eta=\{\eta_{i}:\;M_{i-1}\to M_{i}\to A\dashrightarrow\}_{i=1}^{n}$$ be a family of
$\E$-triangles. Then there exists an $\E$-triangle
$$\xi_{k}:\; M_0\to M_{k}\to A^{k}\dashrightarrow$$
for any $k\in [1,n]:=\{1,2,\cdots,n\}$ for any $n\in\mathbb{Z}^{+}$.
\end{lemma}

\proof We will proceed by induction on $k$. For $k=1,$ we can choose that
$\xi_{1}:=\eta_{1}$ is the required $\E$-triangle.

If $k>1.$ we assume that there exists $\xi_{k-1}$. By (ET4), we get the
following commutative diagram
$$\xymatrix{M_0\ar[r]\ar@{=}[d] & M_{k-1}\ar[r]\ar[d] &A^{k-1}\ar@{-->}[r]\ar[d]&\\
M_0\ar[r] & M_k\ar[r]\ar[d] & L_k\ar@{-->}[r]\ar[d]&\\
& A\ar@{-->}[d]\ar@{=}[r] & A\ar@{-->}[d]\\&&}$$
of $\E$-triangles. Since
$\E(A,A)=0$, the third column $\E$-triangle
of the above diagram splits and then $L_{k}\simeq
A^{k}$. Hence the second row of
the above diagram is the required $\E$-triangle $\xi_{k}$.  \qed
\medskip

Let $(\C,\E,\s) $ be an extriangulated category and $\Theta=\{\Theta(i)\}_{i=1}^{t}$ be a family of
objects in $\C$.  For a given $\Theta$-filtration $\xi=\{\xi_k\colon ~M_{k-1}\to M_k\to X_k\dashrightarrow\}_{k=0}^n$ of $M\in\F(\Theta),$ we denote by
$[M:\Theta(i)]_\xi$ the \emph{$\xi$-filtration multiplicity} of $\Theta(i)$ in $M.$ That is
$[M:\Theta(i)]_\xi$ is the cardinal of the set $\{k\in[0,n]\;|\; X_k\simeq\Theta(i)\}.$
In general, the filtration multiplicity could be depending on a given $\Theta$-filtration. Notice that
$\ell_{\Theta,\xi}(M)=\sum\limits_{k=1}^t\;[M:\Theta(i)]_\xi.$

\begin{lemma}\label{lem7} Let $\Theta=\{\Theta(i)\}_{i=1}^{t}$ be
a family of objects in an extriangulated category $(\C,\E,\s)$, and let
$\leq$ be a linear order on $[1,t]$ such that
$\E(\Theta(j),\Theta(i))=0$ for all $j\geq i$.
If $\xi$ is a $\Theta$-filtration of $M\in \F(\Theta)$, then there exists a $\Theta$-filtration
$\eta$ of $M$ and a family $\Xi$ of $\E$-triangles
satisfying the following conditions.
\begin{enumerate}
\item [\rm (1)] $m(i):=[M:\Theta(i)]_{\xi}=[M:\Theta(i)]_{\eta}$ for all $i\in [1,t].$

\item [\rm (2)] The family $\eta$ is ordered, that is,
$$\eta=\{\eta_{i}\colon~M_{i-1}\to M_{i}\to \Theta(k_{i})\dashrightarrow\}_{i=0}^{n}$$ with
$\Theta(k_{0}):=0$, $M_{-1}:=0$ and $k_{n}\leq k_{n-1}\leq \cdots
\leq k_{1}$.

\item [\rm (3)] $\Xi=\{\Xi_{i}\colon~M_{i-1}'\to M_{i}'\to \Theta(\lambda_{i})^{m(\lambda_{i})}\dashrightarrow\}_{i=0}^{d},$  $\{\Theta(\lambda_{i})\}_{i=1}^{d}$
is the set consisting of the different $\Theta(j)$ appearing in
the $\Theta$-filtration $\xi$ of $M.$

Moreover
$\Theta(\lambda_{0}):=0,$ $M_{-1}':=0,$ $M'_d=M$ and $\lambda_{d}< \lambda_{d-1}<
\cdots < \lambda_{1}$.
\end{enumerate}
\end{lemma}

\proof This is similar to the proof of Proposition 4.7 in \cite{MS}.

Let $\xi$ be a $\Theta$-filtration of $M\in \F(\Theta).$ If $M=0$, this result is trivial. Now we assume
$M\neq 0$.

We begin to show that (1) and (2) by induction on $n:=\ell_{\Theta,\xi}(M).$

If $n=1$, the $\Theta$-filtration $\xi$ is already ordered and hence
$\eta:=\xi$ satisfies the required conditions.

Suppose  $n\geq 2$ and
$\xi:=\{\xi_{i}\colon~M_{i-1}\to M_{i}\to \Theta(k_{i})\dashrightarrow\}_{i=0}^{n}$ be the
$\Theta$-filtration of $M$. Since $\xi':=\xi-\{\xi_{n}\}$ is a
$\Theta$-filtration of $M_{n-1}$  and
$\ell_{\Theta,\xi'}(M_{n-1})=n-1$, by induction there is an ordered
$\Theta$-filtration $\eta'=\{\eta'_{i}\colon~M'_{i-1}\to M'_{i}\to \Theta(k'_{i})\dashrightarrow\}_{i=0}^{n-1}$ of $M_{n-1}$ with
$k'_{n-1}\leq k'_{n-2}\leq \cdots \leq k'_{1}$ and
$[M_{n-1}:\Theta(i)]_{\xi'}=[M_{n-1}:\Theta(i)]_{\eta'}$ for any $i$. If $k_{n}\leq k'_{n-1}$, then $\eta:=\eta'\cup\{\xi_{n}\}$ satisfies the required conditions.

Now we assume that $k'_{n-1}<k_{n}.$  Put $q:=\max\{m\in [1,n-1]\;|\;
k'_{n-m}<k_{n}\}$. Note that the $\Theta$-filtration
$\eta'\cup\{\xi_{n}\}$ is almost the one we want, the only $\E$-triangle
that does not have its ordered multiplicity is precisely the
$\xi_{n}$. This can be rearranged by applying $q$-times
Lemma \ref{lem5} to $\eta'\cup \{\xi_{n}\}$.

(3) In order to construct $\Xi,$ we use the ordered $\Theta$-filtration $\eta$ from (2). Our detailed process is as follows. For any $i\in[1,n]$, we group the $k_{i}$
that are the same and rename them by $\lambda_i$. Therefore we obtain $\lambda_{d}< \lambda_{d-1}< \dots < \lambda_{1}$ and then $\Theta(\lambda_{1}), \cdots, \Theta(\lambda_{d})$ are the different $\Theta(j)$ appearing in the $\Theta$-filtration $\eta$  of $M.$ 

Define
$s(i):=m(\lambda_{i})=[M:\Theta(\lambda_{i})],$ $\alpha(i):=\sum\limits_{j=1}^{i}s(i)$ and $\alpha(0):=-1$.

We divide the filtration $\eta$ into the following pieces
$$\{\eta_{i}\colon~\xymatrix{M_{i-1}\ar[r] &
M_{i}\ar[r] & \theta(\lambda_{q})\dashrightarrow}\}_{i=\alpha(q-1)+1}^{\alpha(q)},$$ where $q\in[1,d]$. By
Lemma \ref{lem6}, for any $q\in[1,d],$ we have the following $\E$-triangle
$$\Xi_{q}\colon~ \xymatrix{M_{\alpha(q-1)}\ar[r] & M_{\alpha(q)}\ar[r] &
\theta(\lambda_{q})^{s(q)}\dashrightarrow}.$$ Hence, by setting  $\Xi_{0}:=\eta_{0}$
and $M_{i}':=M_{\alpha(i-1)}$ for any $i\in[1, d],$ we conclude that the filtration
$\Xi=\{\Xi_{i}\}_{i=0}^{d}$ satisfies the required conditions.  \qed

\medskip
Let $(\C,\E,\s)$ be an extriangulated category and $\Theta=\{\Theta(i)\}_{i=1}^{t}$ be a family of objects in
$\C$.  We denote by $\Theta^\oplus$ the subcategory of $\C$, whose objects are the
finite direct sums of copies of objects in $\Theta$.

\begin{lemma}\label{lem8} Let $(\C,\E,\s)$ be an extriangulated category and $\Theta=\{\Theta(i)\}_{i=1}^{t}$ be a family of objects in
$\C$. Then we have the following statements hold.
\begin{itemize}
\item[\rm (1)] $\F(\Theta)=\F(\Theta^\oplus).$
\item[\rm (2)] If $\C$ is an artin extriangulated $R$-category, then $\Theta^\oplus$ is functorially finite.
\end{itemize}
\end{lemma}

\proof (1) Since $\Theta\subseteq \Theta^{\oplus}$, we have
$\F(\Theta)\subseteq\F(\Theta^{\oplus})$.

Suppose $M\in \F(\Theta^{\oplus})$.  We show that $M\in\F(\Theta)$ by induction on
$m:=\ell_{\Theta^{\oplus}}(M)$.

If $m=1$,
then $M\in \Theta^{\oplus}$ and hence $M=\bigoplus\limits_{i=1}^n\Theta(k_{i})^{m_i}.$
By Lemma \ref{lem2}, we know that $\F(\Theta)$ is closed under finite direct sums. Note that $\Theta(k_{i})\in
\F(\Theta)$,  we get that $M\in\F(\Theta)$.

If $m>1$,  then there exists an $\E$-triangle
$$M_{m-1}\to M\to \Theta(k_{m})^{\lambda(m)}\dashrightarrow$$ where
$\ell_{\Theta^{\oplus}}(M_{m-1})=m-1.$ Thus we obtain that $M_{m-1}\in
\F(\Theta)$ by induction. Therefore $M\in \F(\Theta)$ since
$\F(\Theta)$ is closed under extensions (see Lemma \ref{lem2}).

(2) The proof given in \cite[Proposition 4.2]{AS} can be easily extended to the context of an artin extriangulated $R$-category.

\begin{lemma}\label{lem9} Let $(\C,\E,\s)$ be an extriangulated category, $\X,\Y$ and $\Z$ be three classes of objects in
$\C$. Then
$(\X\ast\Y)\ast\Z=\X\ast(\Y\ast\Z)$. That is to say, the operative $\ast$ is associative.
\end{lemma}

\proof For any object $M\in(\X\ast\Y)\ast\Z$, then there exists an $\E$-triangle
$$L\overset{~}{\longrightarrow}M\overset{~}{\longrightarrow}Z\overset{}{\dashrightarrow}$$
where $L\in\X\ast\Y$ and $Z\in\Z$. Since $L\in\X\ast\Y$, then there exists an $\E$-triangle
$$X\overset{~}{\longrightarrow}L\overset{~}{\longrightarrow}Y\overset{}{\dashrightarrow}$$
where $X\in\X$ and $Y\in\Y$.
By (ET4), we have the following commutative diagram
$$\xymatrix{X\ar[r]\ar@{=}[d]&L\ar[r]\ar[d]&Y\ar[d]\ar@{-->}[r]&\\
X\ar[r]&M\ar[d]\ar[r]&N\ar[d]\ar@{-->}[r]&\\
&Z\ar@{-->}[d]\ar@{=}[r]&Z\ar@{-->}[d]\\
&&}$$
of $\E$-triangles. It follows that $M\in\X\ast(\Y\ast\Z)$ since $X\in\X$ and $N\in\Y\ast\Z$.
This shows that
$(\X\ast\Y)\ast\Z\subseteq\X\ast(\Y\ast\Z)$. Similarly, one can show that $(\X\ast\Y)\ast\Z\supseteq\X\ast(\Y\ast\Z)$.
Hence $(\X\ast\Y)\ast\Z=\X\ast(\Y\ast\Z)$. \qed

\begin{lemma}\label{lem10} Let $\X$ be a class of objects in an extriangulated category
$\C$ such that $0\in \X$ and $\X$ is closed under
isomorphisms. Then $\F_{n}(\X)=\ast_{i=1}^n\,\X$ for any $n\geq 1,$ and
$\F_{k}(\X)\subseteq \F_{k+1}(\X)$ for any $k\in\mathbb{N}.$
\end{lemma}

\proof Note that $\F_{0}(\X):=\{0\}$ and $\X\subseteq
\F_{1}(\X)=\{0\}\ast\X$. Since
$\X$ is closed under isomorphisms, then
$\{0\}\ast\X\subseteq \X$. Hence,
$\F_{1}(\X)=\X$ and so $\F_{2}(\X)=\X\ast\X$. Continuing in the same
way, we get that $\F_{n}(\X)=\ast_{i=1}^n\,\X$ for any $n\geq 1.$ By Lemma \ref{lem9}, we have
$\F_{k+1}(\X)=\X\ast\F_{k}(\X).$ Since $0\in \X,$ we conclude that
$\F_{k}(\X)\subseteq \F_{k+1}(\X)$ for any $k\in\mathbb{N}.$ \qed

\medskip
The following result was proved by Ringel \cite[Theorem 1]{R} in case $\C$ is a module category,
and proved by Mendoza-Santiago in \cite[Theorem 4.10]{MS} in case $\C$ is a triangulated category.
The proof we give here uses the extriangulated version of Gentle-Todorov's theorem due to He
\cite[Theorem 3.3]{H}. More precisely, let $(\C,\E,\s)$ be an extriangulated category with enough projectives,  $\X$ and $\Y$ be two covariantly finite subcategories of $\C$. Then $\X\ast\Y$ is a covariantly finite subcategory of $\C$.

\begin{theorem}\label{main1}
Let $(\C,\E,\s)$ be an artin extriangulated $R$-category with enough projectives and enough injectives,
and $\Theta=\{\Theta(i)\}_{i=1}^{t}$ be
a family of objects in $\C$, and let
$\leq$ be a linear order on the set $[1,t]$ such that
$\E(\Theta(j),\Theta(i))=0$ for all $j\geq i$. Then
$\F(\Theta)=\ast_{i=1}^t\,\Theta^{\oplus}$ and it is functorially finite.
\end{theorem}

\proof Put $\X:=\Theta^{\oplus}.$  We claim that $\F_{t}(\X)=\F(\Theta)$.
In fact, by Lemma \ref{lem8}, we have that
$\F(\X)=\F(\Theta)$ and then
$\F(\Theta)=\bigcup\limits_{n\in
\mathbb{N}}\F_{n}(\X)$.

Suppose $M\in \F(\Theta),$ and we consider a $\Theta$-filtration $\xi$ of $M.$ By Lemma \ref{lem7}(3), there exists a family of $\E$-triangles
$$\Xi=\{\Xi_{i}\colon\xymatrix{M_{i-1}'\ar[r] & M_{i}'\ar[r]
& \Theta(\lambda_{i})^{m(\lambda_{i})}\ar@{-->}[r] & }\}_{i=0}^{d},$$ where $\{\Theta(\lambda_{i})\}_{i=1}^{d}$
is the set of the different $\Theta(j)$ appearing in the $\Theta$-filtration
$\xi$ of $M$, $M_{d}'=M$ and $\lambda_{d}<
\lambda_{d-1}< \cdots < \lambda_{1}$. Therefore $M\in
\F_{d}(\X)$ with $d \leq t$. Since
$\X$ is closed under isomorphisms and contain the zero object,  by Lemma \ref{lem10},
it follows that $\F_{d}(\X)\subseteq
\F_{t}(\X)$. Thus
$\F(\Theta)\subseteq \F_{t}(\X)$ and then
$\F_{t}(\X)=\F(\Theta)$ as required.

By Lemma \ref{lem8}(2), we have that $\X$ is functorially finite. Moreover, from Lemma \ref{lem10} and the claim above, it follows that $\F(\Theta)=\ast_{i=1}^t\,\X.$ Hence this result follows from \cite[Theorem 3.3]{H} and its dual.  \qed
\medskip

As a special case of Theorem \ref{main1}, we have the following.
\begin{corollary}\emph{\cite[Theorem 4.10]{MS}}
Let $\C$ be an artin triangulated $R$-category with shift functor $[1]$
and $\Theta=\{\Theta(i)\}_{i=1}^{t}$ be
a family of objects in $\C$, and let
$\leq$ be a linear order on the set $[1,t]$ such that
${\rm Hom}_{\C}(\Theta(j),\Theta(i)[1])=0$ for all $j\geq i$. Then
$\F(\Theta)=\ast_{i=1}^t\,\Theta^{\oplus}$ and it is functorially finite.
\end{corollary}

\proof Since any triangulated category can be viewed as an extriangulated category with enough projectives and enough injectives.  \qed

\begin{corollary}\emph{\cite[Theorem 1]{R}}
Let $\Lambda$ be an artin algebra and ${\rm mod}(\Lambda)$ the category of finitely generated left
$\Lambda$-modules, and $\Theta=\{\Theta(i)\}_{i=1}^{t}$ be
a family of objects in $\C$, and let
$\leq$ be a linear order on the set $[1,t]$ such that
${\rm Ext}_{\Lambda}^1(\Theta(j),\Theta(i))=0$ for all $j\geq i$. Then
$\F(\Theta)$ is functorially finite.
\end{corollary}

\proof Since any module category can be viewed as an extriangulated category with enough projectives and enough injectives.  \qed

\medskip

Next we introduce the notions of $\Theta$-projective objects and $\Theta$-injective objects.

\begin{definition} \label{def1}
Let $(\C,\E,\s)$ be an extriangulated category and $\Theta=\{\Theta(i)\}_{i=1}^{t}$ be a family of objects in
$\C$. The {\bf $\Theta$-projective} objects in $\C$ is the class $\P(\Theta):={}^{_{\E}\perp}\F(\Theta).$ Dually, the {\bf $\Theta$-injective} objects in $\C$ is the
class $\I(\Theta):=\F(\Theta)^{_{\E}\perp}.$

Note that by Lemma \ref{lem4} and its dual, we have that $\P(\Theta)={}^{_{\E}\perp}\Theta$  and $\I(\Theta)=\Theta^{_{\E}\perp}.$
\end{definition}

\begin{lemma}\label{lem15} Let $(\C,\E,\s)$ be an artin extriangulated $R$-category  and $\Theta=\{\Theta(i)\}_{i=1}^{n}$ be a
family of objects in $\C$ such
that  $\E(\Theta(j),\Theta(i))=0$ for any $j\geq i$.
Consider $t\in [1,n]$ and $N\in \C$  such that
$\E(\Theta(j), N)=0$ for $j>t$. Then there
exists an $\E$-triangle in $\C$
$$\xymatrix{N\ar[r] & N_{t}\ar[r] & \Theta(t)^{m}\ar@{-->}[r] &}$$
where $m:=\ell_{R}(\E(\Theta(t),N))$ and
$\E(\Theta(j), N_{t})=0$ for any $j\geq t$.
\end{lemma}

\proof If $\E(\Theta(t),N)=0$, then the $\E$-triangle
$N\xrightarrow{~1_{N}~}N\xrightarrow{~}0\dashrightarrow$ satisfies the required conditions.

Now we assume $\E(\Theta(t),N)\neq 0.$ By the dual of Lemma \ref{y1}, there exists an $\E$-triangle
$$\eta\colon\xymatrix{N\ar[r] & N_{t}\ar[r] & \Theta(t)^{m}\ar@{-->}[r] &}$$ such that the morphism
$\delta_{\sharp}^{\Theta(t)}\colon\Hom_{\C}(\Theta(t),\Theta(t)^{m})
\to\E(\Theta(t),N)$ is surjective. Applying the functor
$\Hom_{\C}(\Theta(j),-)$  to the $\E$-triangle $\eta,$ we have the following
exact sequence
$$\xymatrix{\Hom_{\C}(\Theta(j),\Theta(t)^{m})\ar[r] & \E(\Theta(j),N)\ar[r] &
\E(\Theta(j),N_{t})\ar[r] &
\E(\Theta(j),\Theta(t)^{m}).}$$ Since
$\E(\Theta(j),\Theta(t))=0$ for any $j\geq t$ and
$\E(\Theta(j), N)=0$ for $j>t$, then
$\E(\Theta(j), N_{t})=0$ for $j>t$. For $j=t$, we know that
$\delta_{\sharp}^{\Theta(t)}$ is an epimorphism and hence
$\E(\Theta(t), N_{t})=0$.  \qed

\begin{lemma}\label{lem16}Let $(\C,\E,\s)$ be an artin extriangulated $R$-category  and $\Theta=\{\Theta(i)\}_{i=1}^{n}$ be a
family of objects in $\C$ such
that  $\E(\Theta(j),\Theta(i))=0$ for any $j\geq i$.
Consider $t\in [1,n]$ and $N\in \C$ such that
$\E(\Theta(j), N)=0$ for $j>t$. Then there exists
an $\E$-triangle in $\C$
$$\xymatrix{N\ar[r] & Y\ar[r] & X\ar@{-->}[r] &}$$ where $X\in
\F(\{\Theta(i)\;|\; i\in [1,t]\})$ and $Y\in\I(\Theta).$
\end{lemma}

\proof Since $\E(\Theta(j),N)=0$ for $j>t$, by Lemma \ref{lem15},
there exists an $\E$-triangle
$$\eta_{t+1}\colon\xymatrix{N\ar[r]^{\mu_{t}} & N_{t}\ar[r] & Q_{t}\ar@{-->}[r] &}$$
with $Q_{t}:=\Theta(t)^{m_{t}}$ and $\E(\Theta(j),N_{t})=0$ for $j\geq t$.
Similarly, there exists an $\E$-triangle
$$\eta_{t}\colon\xymatrix@C=1.2cm{N_{t}\ar[r]^{\mu_{t-1}\quad} & N_{t-1}\ar[r] & Q_{t-1}\ar@{-->}[r] &}$$
with $Q_{t-1}:=\Theta(t-1)^{m_{t-1}}$ and $\E(\Theta(j),N_{t-1})=0$ for $j\geq t-1$. Continue this process, we
obtain $\E$-triangles
$$\eta_{i}\colon\xymatrix{N_{i}\ar[r]^{\mu_{i-1}} & N_{i-1}\ar[r] & Q_{i-1}\ar@{-->}[r] &}$$
where $Q_{i-1}:=\Theta(i-1)^{m_{i-1}}$ and $\E(\Theta(j), N_{i})=0$ for any $j\geq i$. In what follows, for
$\alpha_{r}:=\mu_{t-r}\dots\mu_{t-1}\mu_{t}$ where $0\leq r\leq
t-1$, we will inductively construct $\E$-triangles
$$\xi_{r}\colon\xymatrix{N\ar[r]^{\alpha_{r}\;\;} & N_{t-r}\ar[r] &
X_{t-r}\ar@{-->}[r]&}$$ with $X_{t-r}\in
\F(\{\Theta(i)\;|\; i\in [t-r,t]\})$ and $\E(\Theta(j),N_{t-r})=0$ for
$j\geq t-r$. 

If $r=0$, we put $\xi_{0}:=\eta_{t+1}$. 

Now we assume that $r>0$ and the $\E$-triangle $\xi_{r}$ is already constructed. By (ET4), we have the following commutative diagram
$$\xymatrix{N\ar[r]^{\alpha_{r}}\ar@{=}[d] & N_{t-r}\ar[r]\ar[d]^{\mu_{t-r-1}}
& X_{t-r}\ar@{-->}[r]\ar[d]
&\\
N\ar[r]^{\alpha_{r+1}} & N_{t-r-1}\ar[r]\ar[d] &
X_{t-r-1}\ar@{-->}[r]\ar[d] &  \\
& \Theta(t-r-1)^{m_{t-r-1}}\ar@{=}[r]\ar@{-->}[d] &
\Theta(t-r-1)^{m_{t-r-1}}\ar@{-->}[d]\\
&&}$$
of $\E$-triangles.
By induction, we have that
$X_{t-r}\in \F(\{\Theta(i)\;|\;i\in [t-r,t]\})$. Hence
$$\Theta(t-r-1)^{m_{t-r-1}},~X_{t-r}\in \F(\{\Theta(i)~|~
i\in [t-r-1,t]\}).$$ Since $\F(\{\Theta(i)~|~ i\in
[t-r-1,t]\})$ is closed under extensions, then
$$X_{t-r-1}\in \F(\{\Theta(i)\;|\; i\in [t-r-1,t]\}).$$Furthermore
$\E(\Theta(j),N_{t-r-1})=0$ for any $j\geq t-r-1$.
Thus $\xi_{r+1}$ is the $\E$-triangle from the second row of the
above diagram. This shows that $\xi_{t-1}$ is the required $\E$-triangle.  \qed

\begin{theorem}\label{main2} Let $(\C,\E,\s)$ be an artin extriangulated $R$-category  and $\Theta=\{\Theta(i)\}_{i=1}^{n}$ be a
family of objects in $\C$ such
that  $\E(\Theta(j),\Theta(i))=0$ for any $j\geq i$, and let
$\leq$ be a linear order on the set $[1,t]$ such that
$\E(\Theta(j),\Theta(i))=0$ for all $j\geq i$. Then we have the following statements hold.
\begin{itemize}
 \item[\rm (1)] For any object $X\in\C$, there are two $\E$-triangles in $\C$
$$\xymatrix{X\ar[r] & Y_X\ar[r] & C_X\ar@{-->}[r] & ~~\text{ with } Y_X\in\I(\Theta)~\textrm{and}~~C_X\in\F(\Theta)},$$
$$\xymatrix{K_X\ar[r] & Q_X\ar[r] & X\ar@{-->}[r] & ~~\text{ with } Q_X\in\P(\Theta) ~\textrm{and}~~K_X\in\F(\Theta)}.$$
\item[\rm (2)] $\P(\Theta)$ is a precovering class and $\I(\Theta)$ is a preenveloping one in $\C.$
\end{itemize}
\end{theorem}

\proof (1) For simplicity, we assume that the linear order $\leq$ on the set $[1,t]$ is the natural one. 
Moreover, we only prove that there exists the first $\E$-triangle because the existence of other $\E$-triangle follows by duality.

Let $X\in\C$ and $t:=n.$ By Lemma \ref{lem16}, we have an $\E$-triangle
 $$X\to Y_X\to C_X\dashrightarrow$$ in $\C$ such that $Y_X\in\I(\Theta)$ and $C_X\in\F(\Theta).$
\medskip

(2) Now we show that  $\I(\Theta)$ is a preenveloping class in $\C.$
In fact, for any $X\in \C$, by (1), there exists an
$\E$-triangle
$$\xymatrix{X\ar[r]^{\beta} & Y_X\ar[r] & C_X\ar@{-->}[r] &}$$ with $Y_X\in\I(\Theta)$ and $C_X\in\F(\Theta).$ We claim
that $\beta$ is an $\I(\Theta)$-preenvelope of $X$. Indeed, we consider a morphism $\beta':X \to  Y'$ with $Y'\in\I(\Theta).$ By Lemma \ref{y2}, we have the following
commutative diagram
$$\xymatrix@C=1.2cm{ X\ar[d]^{\beta'}\ar[r]^{\beta} & Y_X\ar@{-->}[dl]_{w}\ar[r]\ar[d]^{\gamma} &
C_X\ar@{=}[d]\ar@{-->}[r]\ar[dl]_{v'}&\\
 Y'\ar[r]^{u} & L\ar[r]^{v} & C_X\ar@{-->}[r]^{\phi}&}$$
of $\E$-triangles. Since  $Y'\in\I(\Theta)=\F(\Theta)^{_{\E}\perp}$ and $C_X\in \F(\Theta),$ we get that
$\E(C_X,Y')=0.$  Therefore $\phi=0$ and then $v$ is a split epimorphism.
That is to say, there exists a morphism $v'\colon C_X\to L$ such that
$vv'=1_{C_X}$.  By Lemma \ref{y0},
we have that $\beta'$ factors through $\beta$.
This shows that $\beta$ is an $\I(\Theta)$-preenvelope of $X.$

By duality, we can show that $\P(\Theta)$ is a precovering class in $\C$ by
using the second $\E$-triangle in (1).  \qed

\begin{remark}
In Theorem \ref{main2}, if $(\C,\E,\s)$ is an artin triangulated $R$-category, Theorem \ref{main2} is just the
Theorem 4.14 in \cite{MS}.
\end{remark}

\textbf{Panyue Zhou}\\
D\'{e}partement de Math\'{e}matiques, Universit\'{e} de Sherbrooke, Sherbrooke,
Qu\'{e}bec J1K 2R1, Canada.\\
and\\
College of Mathematics, Hunan Institute of Science and Technology, 414006, Yueyang, Hunan, P. R. China.\\
E-mail: \textsf{panyuezhou@163.com}


\begin{thebibliography}{99}
\bibitem[AR]{AR} M. Auslander, I. Reiten.
Applications of contravariantly finite subcategories.
Adv. Math. 86(1): 111-152, 1991.


\bibitem[AS]{AS} M. Auslander, S. O. Smal{\o}. Preprojective modules over Artin algebras. J. Algebra 66(1): 61-122, 1980.


\bibitem[CYZ]{CYZ} X. Chen, Y. Ye, P. Zhang.
Algebras of derived dimension zero.
Comm. Algebra 36(1): 1-10, 2008.


\bibitem[H]{H} J. He. Extensions of covariantly finite subcategories revisited. Czechoslovak Math. J. 69(144): 403-415, 2019.

\bibitem[HZZ]{HZZ} J. Hu, D. Zhang, P. Zhou. Proper classes and Gorensteinness in extriangulated categories,
arXiv:1906.10989, 2019.


\bibitem[LN]{LN} Y. Liu, H. Nakaoka.
Hearts of twin cotorsion pairs on extriangulated categories.
J. Algebra, 528: 96-149, 2019.

\bibitem[MS]{MS} O. Mendoza, V. Santiago.
Homological systems in triangulated categories.
Appl. Categ. Structures 24(1): 1-35, 2016.

\bibitem[NP]{NP} H. Nakaoka, Y. Palu.  Extriangulated categories, Hovey twin cotorsion pairs and model structures. Cah. Topol. G¨¦om. Diff¨¦r. Cat¨¦g. 60(2): 117-193, 2019.


\bibitem[R]{R} C.M. Ringel. The category of modules with good filtrations over a quasi-hereditary algebra has almost split sequences. Math. Z. 208(2): 209-233, 1991.

\bibitem[ZZ]{ZZ} P. Zhou, B. Zhu, Triangulated quotient categories revisited. J. Algebra. 502: 196-232, 2018.

\end{thebibliography}
\end{document}